\documentclass[12pt]{amsart}

\usepackage{amsmath,amssymb,amscd,amsthm}

\title[On the fundamental group of an arrangement]%
{ On the fundamental group of the complement of a complex hyperplane
arrangement }

\author{Grigory Rybnikov}
\thanks{This work was supported in part by RFBR grant no. 09-01-12185-ofi\_m, by HSE grant no. 09-09-0010, and by HSE project no. TZ-62.0 ``Mathematical investigations in small-dimensional topology, algebraic geometry, and representation theory.''}

\address{Higher School of Economics, Department of Mathematics}
\email{gr@mccme.ru, grybnikov@hse.ru}
\dedicatory{To the blessed memory of I.~M.~Gelfand}
\date{}

\allowdisplaybreaks
\newcommand{\les}{\leqslant}

\newcommand{\ep}{\varepsilon}
\newcommand{\clll}{\mathcal{L}}
\newcommand{\ppp}{\mathcal{P}}
\newcommand{\rrr}{\mathcal{R}}
\newcommand{\aaa}{\mathcal{A}}
\newcommand{\zz}{\mathbb{Z}}
\newcommand{\qq}{\mathbb{Q}}
\newcommand{\cc}{\mathbb{C}}
\newcommand{\ff}{\mathbb{F}}

\newcommand{\po}{\pi_1}
\newcommand{\g}[1]{\gamma_{#1}}

\DeclareMathOperator{\gr}{gr}

\DeclareMathOperator{\Ker}{Ker}
\DeclareMathOperator{\im}{Im}
\DeclareMathOperator{\Hom}{Hom}
\DeclareMathOperator{\Aut}{Aut}
\newtheorem{theorem}{Theorem}[section]
\newtheorem{lemma}[theorem]{Lemma}
\newtheorem{proposition}[theorem]{Proposition}
\theoremstyle{definition}
\newtheorem{definition}{Definition}
\theoremstyle{remark}
\newtheorem{remark}{Remark}

\newcommand{\prf}[1]{\begin{proof}{#1}\end{proof}}

\begin{document}

\begin{abstract}
We construct two combinatorially equivalent
line arrangements in the complex projective plane
such that the fundamental groups
of their complements are not isomorphic. The proof
uses a new invariant of the fundamental group of the complement to
a line arrangement of a given combinatorial type
with respect to isomorphisms inducing the canonical isomorphism of the first
homology groups.
\end{abstract}

\maketitle

This paper is a revised version of the preprint~\cite{Ry}, which appeared in
1994 and was subsequently cited by several authors.
The invariant used there was later discussed, in particular,
in~\cite{Art}. I~am grateful to V.~A.~Vassiliev, who encouraged me to submit the paper for publication, and to S.~A.~Yuzvisky for advice on the exposition.

I am grateful to the referees for remarks that helped me to correct several shortcomings of the first version of the paper.

\section{Introduction}

Many topological properties of a complex hyperplane arrangement can be
expressed in terms of its combinatorial structure (matroid). For
example, this is so for the cohomology ring of the complement of a complex hyperplane 
arrangement \cite{OS} and for the Malcev completion of its fundamental
group \cite{Ko}.
However, all known algorithms for computing the fundamental group itself
use non-matroidal data. (In the case of a complexified real arrangement, the answer can be expressed
in terms of the corresponding oriented matroid, which requires more information then the mere
dimensions of all intersections of hyperplanes, as for the ordinary matroid.)

Due to a classical theorem of Zariski, the fundamental group of the
complement of a hypersurface in $\cc P^n$ is isomorphic to the
fundamental group of the complement of the hypersurface section in a generic 2-dimensional plane.
Hence, in the case under consideration, it suffices to study the fundamental groups of the complements of line arrangements in
$\cc P^2$. In this paper we describe the combinatorial structure of
line arrangements in terms of incidence structures, which is
essentially equivalent to using the matroid language but has the advantage
of being more visual.

The aim of this paper is to construct two combinatorially equivalent
line arrangements in $\cc P^2$ whose complements have non-isomorphic
fundamental groups. Our approach is as follows.

Let $X$ be the complement of a line arrangement $L$, and let $C$ be the incidence structure describing the combinatorics of~$L$. By $H$ we denote a free Abelian group whose basis is indexed by the lines of~$C$. The first homology group $H_1(X,\zz)$ is a free Abelian roup as well, and its basis consists of cycles going around the lines of~$L$ in the positive direction. Thus, it is canonically isomorphic to~$H$. Let $X'$ be the complement of another line arrangement whose combinatorics is described by the same structure $C$. In view of the canonical isomorphisms $H_1(X,\zz)\cong H$ and $H_1(X',\zz)\cong H$, we have canonical projections $\po(X)\to H$ and $\po(X')\to H$, whose kernels are the commutator subgroups of $\po(X)$ and $\po(X')$.

First, we obtain a necessary condition for the existence of an isomorphism $\po(X)\to\po(X')$ concordant with these projections, i.e., an isomorphism inducing the identity automorphism of~$H$. To this end, for a given structure $C$, we construct an invariant, which depends on the group $\po(X)$ together with the projection $\po(X)\to H$. If this invariant takes different values for $\po(X)$ and $\po(X')$, then no isomorphism $\po(X)\to\po(X')$ of the above type can exist. The invariant depends only on
$\po(X)/\gamma_4\po(X)$, where $\gamma_k G$ is the $k$th term of the lower central series of~$G$.

Next, we consider two complex conjugate realizations of the incidence structure~$C_8$ corresponding to the MacLane matroid \cite{McL} (see also \cite{Z}) and show that our invariant distinguishes between these realizations. We also prove that, in the case of~$C_8$, the group of those automorphisms of~$H$ which are reductions of automorphisms of the group~$\po(X)$ is isomorphic to the group of automorphisms of the incidence structure itself.
Then, at last, we combine these results to prove that the incidence structure
$C_{13}$, which is obtained by gluing together two MacLane structures, has at
least two realizations whose complements have non-isomorphic fundamental groups.

\bigskip

I am deeply grateful for helpful discussions and advice to T.~V.~Alekseyevskaya,
A.~W.~M.~Dress, A.~B.~Goncharov, S.~M.~Lvovsky, N.~E.~Mnev, G.~A.~Noskov, S.~Yu.~Orevkov,
D.~A.~Stone, G.~M.~Ziegler, and, especially, 
I.~M.~Gelfand, without whom this paper would never have been written.

\section{Construction of the Distinguishing Invariant}

\begin{definition} An \emph{incidence structure} is a triple
$C=(\clll,\ppp,\succ)$, consisting of a set of lines $\clll$, a set of points $\ppp$, and a binary relation
$\succ$ (incidence) between $\clll$ and $\ppp$ which satisfies the
following axioms.

(1) For any two distinct lines $l$ and $l'$, there is a unique point
$p$ such that $l\succ p$ and $l'\succ p$.

(2) Any point is incident to at least two lines.
\end{definition}

If lines $l$ and $l'$ and a point $p$ are as in axiom (1), we write
$p=l\cap l'$ and call $p$ the intersection point of the lines
$l$ and $l'$. We say that an incidence structure $C=(\clll,\ppp,\succ)$ is
\emph{non-degenerate} if $|\ppp|>1$, i.e., if not all of the lines meet at the same
point.\footnote{The notion of a non-degenerate incidence structure
is essentially equivalent to the notion of a simple rank 3 matroid
on the set of lines.}

\begin{definition}
A \emph{complex realization} of an incidence structure
$C=(\clll,$ $\ppp, \succ)$ is an injective map $\phi$ from $\clll$
to the set of lines in $\cc P^2$ and from $\ppp$ to the set of
points in $\cc P^2$ such that $\phi p\in\phi l$ if and only if $p\prec l$.
\end{definition}

Let $C=(\clll,\ppp,\succ)$ be a finite non-degenerate incidence structure. We order the elements of $\clll$ as $\clll=\{l_0,l_1,\dots,l_n\}$ and call $l_0$ the line at infinity. Let $\ppp_0=\{p\in\ppp\mid p\not\prec l_0\}$. Below we describe in abstract terms the presentation by generators and relations for the fundamental group of the complement of a line arrangement in $\cc^2$ as proposed in \cite{A} (see also \cite{OT}).

Let $F$ be a free group with free generators $w_1,\dots,w_n$, and let $\aaa=\{(i,p)\in\{1,\dots,n\}\times \ppp_0\mid p\prec l_i\}$. Suppose that, for any pair $(i,p)\in\aaa$, we are given an element $w_i(p)\in F$ conjugate in $F$ to $w_i$, i.e., such that $w_i(p)={(g(i,p))}^{-1}w_ig(i,p)$ for some $g(i,p)\in F$. Thus, the elements $w_i(p)$ are determined by some mapping $g:\aaa\to F$.

Let $p\in \ppp_0$ and let $\{i_1,\dots,i_k\}=\{i\in\{1,\dots,n\}\mid l_i\succ p\}$,
where $i_1<\dots<i_k$. We set $c_s(p)=w_{i_s}(p)\cdot
w_{i_{s-1}}(p)\cdot\dots\cdot w_{i_1}(p)\cdot
w_{i_k}(p)\cdot\dots\cdot w_{i_{s+1}}(p)$ for $s=1,\dots,k$
and $r(i_s,p)={c_{s-1}(p)}^{-1}\cdot c_s(p)$
for $s=1,\dots,k$, where $c_0(p)=c_k(p)$. It is readily seen that
$r(i_s,p)=[w_{i_s},c_s(p)]$ and $\prod_{s=1}^k r(i_s,p)=1$,
hence only $k-1$ of the elements $r(i_s,p)$ are independent. For each
$p\in \ppp_0$, consider the set $\{r(i_2,p),\dots,r(i_k,p)\}$ of independent relators; let $\rrr=\rrr(g)=\{r(i,p)\mid l_i\succ p,\,i\neq\min\{j\mid l_j\succ p\}\,\}$ be the union of these sets.
We set $G=G(\rrr)=F/\langle\rrr\rangle$, where $\langle\rrr\rangle$
is the minimal normal subgroup of~$F$ containing~$\rrr$.

Let $\phi$ be a complex realization of the incidence structure $C$, and
let $X=\cc P^2\setminus\bigcup_{l\in\clll}\phi l$. Then, according to~\cite{A},
$\po(X)\cong G(\rrr(g))$ for some particular choice of $g$. Under this isomorphism, the image of $w_i$ corresponds to a small loop passing around $\phi l_i$ in positive direction (the lines in $\cc^2$ have canonical
coorientation).

Consider the lower central series $(\g kF)$ of the group
$F$ ($\g1F=F$, $\g{k+1}F=[F,\g kF]$). It is well known that the groups
$\gr_kF=\g kF/\g{k+1}F$ are free Abelian groups of finite rank.
Moreover, $\gr F=\bigoplus_{k=1}^\infty\gr_kF$ is naturally isomorphic
(as a graded Abelian group) to the free Lie algebra $L$ over $\zz$
with $n$ generators $x_1,\dots,x_n$ which are the images of $w_1,\dots,w_n$ in $F/\g2F$; $L=\bigoplus_{k=1}^\infty L_k$
is graded by the degree of monomials (see, e.g., \cite{MKS}, Chap.~5). Under this isomorphism, the
Lie algebra commutator in $L$ corresponds to the group commutator
in $F$. Hereafter we identify $\gr_kF$ with $L_k$. Let $H=L_1=F/\g{2}F$;
this is a free Abelian group with basis $x_1,\dots,x_n$.
If $G=\po(X)$ for some complex realization of~$C$, then we have
$H\cong H_1(X,\zz)$, and the basis $(x_i)$ is dual to the basis of~$H^1(X,\zz)$ formed by the standard generators of the Orlik--Solomon ring~\cite{OS}.

The group $\langle\rrr\rangle$ is generated by the elements $r(i,p)$ (which are commutators) and by all elements of the form $[w_{i_1}^{\pm 1},[w_{i_2}^{\pm 1}, \dots,[w_{i_k}^{\pm 1}, r(i,P)]\dots]]$ (where the brackets denote the group commutator). Therefore, $\langle\rrr\rangle\subset\g2F$. It is readily seen that the images of the elements of $\rrr$ in $L_2$ are $\bar{r}(i,p)=[x_i,\sum_{j:p\prec l_j}x_j]$ (where the brackets denote the Lie commutator); thus, they depend only on $C$.

In what follows, we need a broader class of groups~$G$. Let $\rrr=\{r(i,p)\mid(i,p)\in\aaa,\, i\neq\min\{j\mid l_j\succ p\}\,\}$  be a set consisting of arbitrary elements of $\g2F$ with the
only condition that the images of $r(i,p)$ in $L_2$ are the same $\bar{r}(i,p)$ as above. We say that such a set $\rrr$ is \emph{admissible} for~$C$. As above, we set $G=G(\rrr)=F/\langle\rrr\rangle$. Since $\langle\rrr\rangle\subset\g2F$, the canonical projection
$F\to H$ is factored through the uniquely determined surjective homomorphism $G\to H$.
The kernel of this homomorphism is the commutator subgroup of~$G$.
We are interested in necessary conditions under which two groups of the form $G(\rrr)$ (determined by different admissible sets $\rrr$) are isomorphic in such a way that the isomorphism is concordant with the canonical projections $G(\rrr)\to H$ specified above.

The lower central series of~$F$ determines the filtration $(\langle\rrr\rangle\cap\g{k}F)$ of the subgroup $\langle\rrr\rangle\subset F$; the corresponding filtration on the quotient group
$G=F/\langle\rrr\rangle$ coincides with its lower central series $(\g kG)$. By
$\gr_k\langle\rrr\rangle$ and $\gr_kG$ denote the graded quotients
of these filtrations. By the definition of the admissible set $\rrr$, the Abelian group
$R_2=\gr_2\langle\rrr\rangle\subseteq L_2$ is generated by the elements
$\bar{r}(i,p)=[x_i,\sum_{j:p\prec l_j}x_j]$; thus, it is uniquely determined by $C$.
Let us show that $R_3=\gr_3\langle\rrr\rangle\subseteq L_3$ is also uniquely determined by $C$.

\begin{proposition}
For any admissible set $\rrr$, the relation $R_3=[H,R_2]$ holds in the Lie algebra $L=\gr F$.
\end{proposition}

\prf{
First, note that the elements $\bar{r}(i,p)$ are linearly independent in $L_2$. Among all elements of the form $[w_j,r(i,P)]$ we choose $u_1,\dots,u_m$ so that their images form a basis in $[H,R_2]$. By analogy with the proof of P.~Hall's basis theorem (see~\cite{MKS}, Theorem~5.13A), it is easy to show that any element of~$\langle\rrr\rangle$ can be uniquely expressed in the form $\prod_{r\in\rrr}r^{c_r}\cdot\prod_{i=1}^mu_i^{d_i}\cdot z$, where $c_r,d_i\in\zz$ and $z\in\langle\rrr\rangle\cap\g{4}F$ (the order of the elements $r\in\rrr$ in the product $\prod_{r\in\rrr}r^{c_r}$ is fixed in advance). Since the images of the elements $r\in\rrr$ in~$\gr_2\langle\rrr\rangle$ are linearly independent, each element of~$\langle\rrr\rangle\cap\g{3}F$
is of the form~$\prod_{i=1}^mu_i^{d_i}\cdot z$. It follows that $R_3=[H,R_2]$.}

We set $P_2=L_2/R_2$ and $P_3=L_3/R_3$. A presentation of a group $G$ by generators and relations as
$G(\rrr)$ for some admissible set $\rrr$ determines an isomorphism $G/\g2G=\gr_1G\cong H$. It is readily seen that, given such an isomorphism, isomorphisms $\g2G/\g3G=\gr_2G\cong P_2$ and $\g3G/\g4G=\gr_3G\cong P_3$ are uniquely determined. In this sense, the groups $\g2G/\g3G$ and $\g3G/\g4G$ are canonically determined by the incidence structure~$C$.

The group $G/\g3G$ is an extension of the group $H$ (which is free Abelian) by the trivial $H$-module $P_2$.
Such extensions are parameterized by the group
$H^2(H,P_2)\cong \Hom(\Lambda^2H,P_2)$. But since $\Lambda^2H\cong L_2$,
we have a natural projection $\chi_2:\Lambda^2H\to P_2$. It is easy
to see that $\chi_2$ is exactly the element of $H^2(H,P_2)$
corresponding to the group extension $1\to P_2\to
G/\g3G\to H\to 1$. Hence $G/\g3G$ is canonically determined by the
incidence structure up to an automorphism trivial on~$H$ and~$P_2$.
The group of such automorphisms is isomorphic to~$\Hom(H,P_2)$.

The group $M=\g2G/\g4G$ is Abelian, but it has non-trivial $H$-module structure. Arguments similar to those used for $G/\g3G$ show that this module is also canonically determined by $C$ up to an action of the Abelian group $\Hom(P_2,P_3)$.

Now, consider the group extension $1\to M\to
G/\g4G\to H\to 1$. It corresponds to a cohomology class $\chi_3\in H^2(H,M)$.
From the short exact sequence of $H$-modules
$0\to P_3\to M\to P_2\to0$
we obtain the long exact sequence
$$
\dots\to H^1(H,P_2)\stackrel{\delta}{\to}H^2(H,P_3)
\stackrel{\alpha}{\to}H^2(H,M)
\stackrel{\beta}{\to}H^2(H,P_2)
\to\dots.
$$
Let us rewrite it in the form
$$
\dots\to\Hom(H,P_2)\stackrel{\delta}{\to}\Hom(\Lambda^2H,P_3)
\stackrel{\alpha}{\to}H^2(H,M)
\stackrel{\beta}{\to}\Hom(\Lambda^2H,P_2)
\to\dots.
$$
We have $\beta(\chi_3)=\chi_2$; thus, $\chi_3$ belongs to $\beta^{-1}(\chi_2)$,
which is a principal homogeneous space over the Abelian group
$\Hom(\Lambda^2H,P_3)/\delta\Hom(H,P_2)$.

Recall that the $H$-module $M$ is determined by $C$ only up to an action of the
Abelian group~$\Hom(P_2,P_3)$; hence only the orbit $\bar{\chi}_3$ of $\chi_3$
under this action has an invariant sense. The corresponding set of orbits $Y=\beta^{-1}(\chi_2)/\Hom(P_2,P_3)$
is a principle homogeneous space over the group $T=\Hom(R_2,P_3)/\bar{\delta}\Hom(H,P_2)$, where $\bar{\delta}$ is the composition of $\delta$ with the natural projection $\Hom(\Lambda^2H,P_3)\to\Hom(R_2,P_3)$.
We see that $Y$, $T$, and the action of $T$ on $Y$ are determined by
$C$ canonically; thus, for any admissible $\rrr$, the element $\bar{\chi}_3(\rrr)\in Y$ is well defined.

Assume now that we have two admissible sets $\rrr$ and $\rrr'$.
Consider the corresponding groups $G=G(\rrr)$ and $G'=G(\rrr')$.
Let $\kappa(\rrr,\rrr')=\bar{\chi}_3(\rrr')-\bar{\chi}_3(\rrr)\in T$.
Clearly, by definition, $\kappa(\rrr,\rrr'')=\kappa(\rrr,\rrr')+\kappa(\rrr',\rrr'')$
for any admissible $\rrr$, $\rrr'$, and $\rrr''$.

The groups $G/\g2G$ and $G'/\g2G'$ are both canonically isomorphic to
$H$. Thus, we have a canonical isomorphism $\lambda_1:G/\g2G\to G'/\g2G'$.
As mentioned above, there always exists an isomorphism $\lambda_2:G/\g3G\to G'/\g3G'$
extending $\lambda_1$, although it is not determined canonically.

\begin{theorem}\label{inv_main}
There exists an isomorphism $\lambda_3:G/\g4G\to G'/\g4G'$ extending
$\lambda_1$ if and only if $\kappa(\rrr,\rrr')=0$
in $T$.
\end{theorem}

\prf{
As shown above, the Lie algebras $\gr(G/\g4G)$ and $\gr(G'/\g4G')$ are canonically isomorphic and are generated by $G/\g2G$ and $G'/\g2G'$, respectively. Therefore, if $\lambda_3$ exists, then $\gr\lambda_3$ is uniquely determined and coincides with this canonical isomorphism. Hence $\kappa(\rrr,\rrr')=0$.

Conversely, if $\kappa(\rrr,\rrr')=0$, then the extensions $1\to M\to G/\g4G\to H\to 1$ and $1\to M\to G'/\g4G' \to H\to 1$ differ by an automorphism of $M$; thus, they become equivalent after applying this automorphism. Hence there exists an isomorphism $\lambda_3$ extending $\lambda_1$. }

In order to use the invariant provided by this theorem, we must describe
$P_2$ and $P_3$. In the general case, they can be described as follows.

Let us identify $L_2$ with $\Lambda^2H$ and denote the basis
of~$H^*$ dual to the basis $(x_i)$ of~$H$ by $(x_i^*)$. Then $L_2^*=\Lambda^2H^*$,
and $R_2^\perp\subset\Lambda^2H^*$ is generated by the elements
$\omega_{ijk}=(x_i^*-x_j^*)\wedge(x_j^*-x_k^*)
=x_i^*\wedge x_j^*+x_j^*\wedge x_k^*+
x_k^*\wedge x_i^*$ for $l_i\succ l_j\cap l_k\in\ppp_0$ and
$\omega_{ij}=x_i^*\wedge x_j^*$
for $l_0\succ l_i\cap l_j$. It is easy to show that $(R_2^\perp)^\perp=
R_2$; hence $P_2$ is a free Abelian group, and $P_2^*\cong R_2^\perp$.
Therefore, we may regard the forms $\omega_{ijk}$ and $\omega_{ij}$
defined above as elements of $P_2^*$.

We have the exact sequence of Abelian groups
$$
0\to \Lambda^3H\stackrel{d}{\to}H\otimes L_2
\stackrel{c}{\to}L_3\to0,
$$
where $d(x\wedge y\wedge z)=x\otimes[y,z]+y\otimes[z,x]+z\otimes[x,y]$
and $c(x\otimes f)=[x,f]$. As we know, $R_3=[H,R_2]$, and hence
$R_3=c(H\otimes R_2)$. Passing to the dual exact sequence
$$
0\leftarrow \Lambda^3H^*\stackrel{\ d^*}{\leftarrow}H^*\otimes\Lambda^2H^*
\stackrel{\ c^*}{\leftarrow}L_3^*\leftarrow0,
$$
we see that $R_3^\perp=\Ker d^*|_{H^*\otimes R_2^\perp}$. Therefore,
$R_3^\perp$ contains all elements of the form
$S_{ijk}=(x_i^*-x_j^*)\otimes\omega_{ijk}$ for
$l_i\succ l_j\cap l_k\in\ppp_0$ and $S_{ij}=x_i^*\otimes\omega_{ij}$
for $l_0\succ l_i\cap l_j$. I do not claim that these elements
generate $R_3^\perp$---as we shall see in the next section, this is not so,---but I cannot give any general description for the remaining generators of $R_3^\perp$.

Translating the definition of $\kappa(\rrr,\rrr')$ from the homological language, in which it is given above, we obtain the following explicit algorithm for computing this invariant.

Let $\rrr=\{r(i,p)\mid(i,p)\in\aaa,\, i\neq\min\{j\mid l_j\succ p\}\,\}$ and
$\rrr'=\{r'(i,p)\mid(i,p)\in\aaa,\, i\neq\min\{j\mid l_j\succ p\}\,\}$
be arbitrary admissible sets. The set $\bar\rrr=\{\bar r(i,p)\mid(i,p)\in\aaa,\, i\neq\min\{j\mid l_j\succ p\}\,\}$ is a basis of $R_2$. To each element $\bar r(i,p)\in\bar\rrr$ we assign the element
$t(i,p)=(r'(i,p))^{-1}r(i,p)\in F$. Since $r(i,p),r'(i,p)\in\g2F$
and their images in $\gr_2F$ coincide (and are equal to $\bar r(i,p)$), we have $t(i,p)\in\g3F$.
We define a homomorphism $R_2\to P_3$ by its action on the basis $\bar\rrr$ as follows: each element
$\bar r(i,p)$ goes to the image of $t(i,p)$ in $P_3=(\gr_3F)/R_3$. The image of this homomorphism in~$T=\Hom(R_2, P_3)/\bar{\delta}\Hom(H,P_2)$ is $\kappa(\rrr,\rrr')$.

We also need the following explicit formula for $\bar\delta$:
$$
\bar\delta f(\bar r(i,p))=\Bigl[f(x_i),\sum_{j:p\prec l_j}x_j\Bigr]+\Bigl[x_i,\sum_{j:p\prec l_j}f(x_j)\Bigr]
$$
for each $f\in\Hom(H,P_2)$.

Let us describe $\kappa(\rrr,\rrr')$ in the case where $\rrr=\rrr(g)$ and $\rrr'=\rrr(g')$ for some mappings
$g,g':\aaa\to F$.

\begin{proposition}
If $g(i,p)\equiv g'(i,p)\pmod{\g2F}$ for all $i=1,\dots,n$ and
$p\prec l_i$, then $\kappa(\rrr,\rrr')=0$. In other words,
$\bar{\chi}_3(\rrr(g))$ depends only on $\bar{g}:\aaa\to H$.
\end{proposition}

\prf{
Indeed, in this case we have $w_i(p)\equiv w'_i(p)\pmod{\g3F}$ and, therefore,
$r(i,p)\equiv r'(i,p)\pmod{\g4F}$. Hence $(r(i,p))^{-1}r'(i,p)\in\g4F$ and $\kappa(\rrr,\rrr')=0$.
}

Consider the Abelian group $A=H^\aaa$. For $a\in A$, we define
$\tilde{\tau}a\in\Hom(R_2,P_3)$ by setting
$$
\tilde{\tau}a(\bar{r}(i,p))=
\Bigl[\left[x_i,a(i,p)\right],\sum_{j:p\prec l_j}x_j\Bigr]+\Bigl[x_i,\sum_{j:p\prec l_j}\left[x_j,a(j,p)\right]\Bigr]
+R_3.
$$
Let $\tau a=\tilde{\tau}a+\bar{\delta}\Hom(H,P_2)$.

\begin{proposition}\label{inv_formula}
If $\rrr=\rrr(g)$ and $\rrr'=\rrr(g')$, then
$\kappa(\rrr,\rrr')=\tau(\bar{g}-\bar{g'})$.
\end{proposition}

\prf{
The required assertion follows from the algorithm for computing $\kappa(\rrr,\rrr')$.
}

It would be interesting to find $\Ker\tau$. Below we give some partial results
in this direction.

Let us define the following elements of $A$:
\ $a^{(0)}_{i,p}$ for $i\in\{1,\dots,n\}$, $p\in\ppp_0$, $p\prec l_i$;
\ $a^{(1)}_{i,p}$ for $i\in\{1,\dots,n\}$, $p\in\ppp_0$;
\ $a^{(2)}_{i,p_1,p_2}$ for $i\in\{1,\dots,n\}$, $p_1,p_2\in\ppp_0$,
$p_1,p_2\prec l_i$ ($p_1$ and $p_2$ are not necessary distinct) by setting
$$
a^{(0)}_{i,p}(j,q)=\delta_{i,j}\delta_{p,q}x_i, \qquad
a^{(1)}_{i,p}(j,q)=\delta_{p,q}x_i, \qquad
a^{(2)}_{i,p_1,p_2}(j,q)=\delta_{i,j}\delta_{p_1,q}\sum_{k:p_2\succ l_k}x_k.
$$
By $U$ denote the subgroup of $A$ generated by them.

\begin{proposition}\label{inv_ker_tau}
The inclusion $U\subseteq\Ker\tilde{\tau}$ holds.
\end{proposition}

\prf{
This inclusion is verified by direct computation.
}

Let $B$ be the subgroup of $A$ consisting of those functions which do 
not depend on $p$.

\begin{proposition}\label{inv_tau_delta}
The inclusion $B\subseteq\tilde{\tau}^{-1}(\bar{\delta}\Hom(H,P_2))$ holds.
\end{proposition}

\prf{
Let $a\in B$, $a(i,p)=\sum_{j=1}^nk_{ij}x_j$. Consider the homomorphism $f\in\Hom(H,P_2)$ defined by
$f(x_i)=\sum_{j=1}^nk_{ij}[x_i,x_j]$. We have $\tilde\tau a = \bar\delta(f)$.
}

Let $W=A/(U+B)$. The homomorphism $\tau$ is the composition of the projection $\pi:A\to W$ and a homomorphism $\bar{\tau}:W\to T$. I do not know whether $\bar{\tau}$ is always injective, but this is so in the case that we consider in the next section.

\section{Computations for the MacLane Incidence Structure}

The MacLane matroid $ML_8$ can be defined as the affine plane over $\ff_3$ minus one element. The lines of the corresponding incidence structure are in one-to-one correspondence with the elements of $\ff_3^2\setminus\{(0,0)\}$, and three lines of $C_8$ meet at one point if and only if the corresponding elements belong to one affine line (in the sense of $\ff_3^2$). It follows that the automorphism group of $C_8$ is isomorphic to $GL_2(\ff_3)$.

Let us index the lines of $C_8$ so that they correspond to the following elements of $\ff_3^2$:
\ $l_0\leftrightarrow(2,2)$, \ $l_1\leftrightarrow(0,1)$, \ $l_2\leftrightarrow(1,0)$,
\ $l_3\leftrightarrow(2,0)$, \ $l_4\leftrightarrow(2,1)$, \ $l_5\leftrightarrow(1,2)$,
\ $l_6\leftrightarrow(0,2)$, \ $l_7\leftrightarrow(1,1)$. Then the points of~$C_8$ are:
 $p_{012}$, $p_{034}$, $p_{056}$, $p_{07}$, $p_{135}$, $p_{147}$, $p_{16}$, $p_{23}$, $p_{246}$,
$p_{257}$, $p_{367}$, $p_{45}$, where $p_{ijk}\prec l_i,l_j,l_k$ and $p_{ij}\prec l_i,l_j$.
Thus $\ppp_0$ consists of the last eight points among those listed above.

It is easy to show that any complex realization of $C_8$ is projectively equivalent
to a realization of the form
\begin{eqnarray*}
\phi l_0 & = & \{(z_0:z_1:z_2)\in\cc P^2\mid z_0=0\}, \\
\phi l_1 & = & \{(z_0:z_1:z_2)\in\cc P^2\mid z_1=0\}, \\
\phi l_2 & = & \{(z_0:z_1:z_2)\in\cc P^2\mid z_1=z_0\}, \\
\phi l_3 & = & \{(z_0:z_1:z_2)\in\cc P^2\mid z_2=0\}, \\
\phi l_4 & = & \{(z_0:z_1:z_2)\in\cc P^2\mid z_2=z_0\}, \\
\phi l_5 & = & \{(z_0:z_1:z_2)\in\cc P^2\mid z_2+\omega z_1=0\}, \\
\phi l_6 & = & \{(z_0:z_1:z_2)\in\cc P^2\mid z_2+\omega z_1=(\omega+1)z_0\}, \\
\phi l_7 & = & \{(z_0:z_1:z_2)\in\cc P^2\mid (\omega+1)z_1+z_2=z_0\},
\end{eqnarray*}
where $\omega$ is a root of the polynomial $x^2+x+1$. The points $\phi p$ for $p\in\ppp$ are defined as the intersection points of the corresponding lines.

We denote the realization of $C_8$ corresponding to
$\omega=\exp(2\pi i/3)$ (to $\omega=\exp(-2\pi i/3)$) by $\phi^+$ (respectively,\ by $\phi^-$). We set $G^+=\po(X^+)/\gamma_4\po(X^+)$ and $G^-=\po(X^-)/\gamma_4\po(X^-)$

\begin{theorem}\label{inv_mclane}
There exists no isomorphism
$G^+\to G^-$
such that the induced isomorphism $H_1(X^+,\zz)\to H_1(X^-,\zz)$ maps
the elements of the canonical basis of $H_1(X^+,\zz)$ to the corresponding
elements of the canonical basis of $H_1(X^-,\zz)$.
\end{theorem}

\prf{
Let $X^\pm=\cc P^2\setminus
\bigcup\phi^\pm l$. Applying Arvola's algorithm, we obtain
the following sets of defining relators for $\po(X^\pm)$
(both in the free group $F$ with generators $w_1,\dots,w_7$):
for $\po(X^+)$
\begin{gather*}
[w_5,w_3,w_1],
\quad
[w_7,w_4,w_1],
\quad
[w_6,w_1],
\quad
[w_3,w_2],
\quad
[w_6,w_4,w_2],
\\
[w_7,w_5,w_5^{-1}w_2w_5],
\quad
[w_7,w_6,w_6^{-1}w_3w_6],
\quad
[w_5,w_3w_6w_4w_6^{-1}w_3^{-1}];
\end{gather*}
and for $\po(X^-)$
\begin{gather*}
[w_5,w_3,w_1],
\quad
[w_7,w_4,w_1],
\quad
[w_6,w_1],
\quad
[w_3,w_5^{-1}w_2w_5],
\quad
[w_7^{-1}w_6w_7,w_4,w_2],
\\
[w_7,w_5,w_5^{-1}w_2w_5],
\quad
[w_7,w_6,w_6^{-1}w_3w_6],
\quad
[w_5,w_7w_4w_7^{-1}].
\end{gather*}
These sets of relators correspond to the mappings $g^\pm:\aaa\to F$ which take the values
\begin{gather*}
g^+(3,p_{367})  =  w_6,
\quad
g^+(4,p_{45})  =  w_6^{-1}w_3^{-1},
\quad
g^+(2,p_{257})  =  w_5,
\quad
g^-(3,p_{367})  =  w_6,
\\
g^-(4,p_{45})  =  w_7^{-1},
\quad
g^-(2,p_{23})  =  w_5,
\quad
g^-(2,p_{257})  =  w_5,
\quad
g^-(6,p_{246})  =  w_7,
\end{gather*}
and the value $1$ at all other pairs.
Let $a_0=\bar{g^+}-\bar{g^-}:\aaa\to H$. The only non-zero values of this mapping are
$$
a_0(4,p_{45}) = x_7-x_6-x_3,
\quad
a_0(2,p_{23}) = -x_5,
\quad
a_0(6,p_{246}) = -x_7.
$$
By Theorem~\ref{inv_main} and Proposition~\ref{inv_formula}, it suffices to show that
$\tau a_0\neq0$ in $T$.

\begin{lemma}
The image of $a_0$ in $W$ is nonzero.
\end{lemma}
\prf{
First, we compute $\tilde{W}=A/U$. Since each generator
$a^{(\alpha)}_{i,p}$ of $U$ vanishes at $(j,q)$ with $q\neq p$,
we have the decompositions $A=\bigoplus_{p\in\ppp_0}A_p$, \ $U=\bigoplus_{p\in\ppp_0}U_p$, and
$\tilde{W}=\bigoplus_{p\in\ppp_0}\tilde{W}_p$,
where $A_p$, $U_p$, and $\tilde{W}_p$ are defined in the obvious way.
For each $p$, the computation is quite simple. All of the $\tilde{W}_p$
turn out to be free Abelian groups; thus, $\tilde{W}$ is free Abelian as well,
and hence its dual $\tilde{W}^*$ is naturally isomorphic to $U^\perp\subset A^*$.

Let $e_{ij}(p)\in A^*$ be defined by $(e_{ij}(p),a)=(x_j^*,a(i,p))$
for $a\in A$, where $(x_i^*)$ is the basis of $H^*$ dual to the basis
$(x_i)$ of $H$. The basis of $U^\perp$ consists of the elements
\begin{eqnarray*}
I(p_{135}) & = & e_{13}(p_{135})-e_{15}(p_{135})+e_{35}(p_{135})-e_{31}(p_{135})+e_{51}(p_{135})-e_{53}(p_{135}),\\
I(p_{147}) & = & e_{14}(p_{147})-e_{17}(p_{147})+e_{47}(p_{147})-e_{41}(p_{147})+e_{71}(p_{147})-e_{74}(p_{147}),\\
I(p_{367}) & = & e_{36}(p_{367})-e_{37}(p_{367})+e_{67}(p_{367})-e_{63}(p_{367})+e_{73}(p_{367})-e_{76}(p_{367}),\\
I(p_{257}) & = & e_{25}(p_{257})-e_{27}(p_{257})+e_{57}(p_{257})-e_{52}(p_{257})+e_{72}(p_{257})-e_{75}(p_{257}),\\
I(p_{246}) & = & e_{24}(p_{246})-e_{26}(p_{246})+e_{46}(p_{246})-e_{42}(p_{246})+e_{62}(p_{246})-e_{64}(p_{246}),\\
J(p_{16}) & = & e_{12}(p_{16})-e_{62}(p_{16})+e_{64}(p_{16})-e_{14}(p_{16})+e_{17}(p_{16})-e_{67}(p_{16})\\
& & {}+e_{63}(p_{16})-e_{13}(p_{16})+e_{15}(p_{16})-e_{65}(p_{16}),\\
J(p_{45}) & = & e_{43}(p_{45})-e_{53}(p_{45})+e_{51}(p_{45})-e_{41}(p_{45})+e_{47}(p_{45})-e_{57}(p_{45})\\
& & {}+e_{52}(p_{45})-e_{42}(p_{45})+e_{46}(p_{45})-e_{56}(p_{45}),\\
J(p_{23}) & = & e_{21}(p_{23})-e_{31}(p_{23})+e_{35}(p_{23})-e_{25}(p_{23})+e_{27}(p_{23})-e_{37}(p_{23})\\
& & {}+e_{36}(p_{23})-e_{26}(p_{23})+e_{24}(p_{23})-e_{34}(p_{23}),\\
K_1(p_{135}) & = & e_{12}(p_{135})+e_{36}(p_{135})-e_{37}(p_{135})+e_{57}(p_{135})-e_{52}(p_{135})-e_{56}(p_{135}),\\
K_2(p_{135}) & = & e_{14}(p_{135})-e_{17}(p_{135})+e_{37}(p_{135})-e_{34}(p_{135})-e_{36}(p_{135})+e_{56}(p_{135}),\\
K_1(p_{147}) & = & e_{12}(p_{147})+e_{13}(p_{147})-e_{15}(p_{147})-e_{43}(p_{147})+e_{75}(p_{147})-e_{72}(p_{147}),\\
K_2(p_{147}) & = & e_{12}(p_{147})+e_{46}(p_{147})-e_{42}(p_{147})-e_{43}(p_{147})+e_{73}(p_{147})-e_{76}(p_{147}),\\
K_1(p_{367}) & = & e_{31}(p_{367})-e_{34}(p_{367})-e_{35}(p_{367})+e_{65}(p_{367})+e_{74}(p_{367})-e_{71}(p_{367}),\\
K_2(p_{367}) & = & e_{34}(p_{367})+e_{62}(p_{367})-e_{64}(p_{367})-e_{65}(p_{367})+e_{75}(p_{367})-e_{72}(p_{367}),\\
K_1(p_{257}) & = & e_{21}(p_{257})+e_{53}(p_{257})-e_{51}(p_{257})-e_{56}(p_{257})+e_{76}(p_{257})-e_{73}(p_{257}),\\
K_2(p_{257}) & = & e_{24}(p_{257})-e_{21}(p_{257})-e_{26}(p_{257})+e_{56}(p_{257})+e_{71}(p_{257})-e_{74}(p_{257}),\\
K_1(p_{246}) & = & e_{21}(p_{246})+e_{47}(p_{246})-e_{41}(p_{246})-e_{43}(p_{246})+e_{63}(p_{246})-e_{67}(p_{246}),\\
K_2(p_{246}) & = & e_{25}(p_{246})-e_{27}(p_{246})+e_{43}(p_{246})+e_{67}(p_{246})-e_{63}(p_{246})-e_{65}(p_{246}).
\end{eqnarray*}

Consider the homomorphism $t:\tilde{W}\to
\zz/3\zz$ defined by
$$
t=2I(p_{367})+J(p_{45})+K_1(p_{135})-K_2(p_{147})+K_1(p_{257})
-K_1(p_{246})+3\zz.
$$
It is not hard to check that $t|_B=0$ and $t(a_0)=1$.
Hence the image of $a_0$ in $W=\tilde{W}/B$ is nonzero. }

\begin{lemma}
The equality $U=\Ker\tilde{\tau}$ holds true for the incidence structure $C_8$.
\end{lemma}

\prf{
To prove the required assertion, we first study $P_3$. It turns out that $R_3^\perp$ has a basis consisting of the elements $S_{ij}=x_i^*\otimes x_i^*\wedge x_j^*$ for $l_0\succ l_i\cap l_j$,
\ $S_{ijk}=(x_i^*-x_j^*)\otimes(x_i^*-x_j^*)\wedge(x_j^*-x_k^*)$ for
$l_i\succ l_j\cap l_k\in\ppp_0$ with $i<j$, $i<k$, and the additional
elements
\begin{eqnarray*}
T_0 & = &
x_7^*\otimes(x_1^*-x_3^*)\wedge(x_3^*-x_5^*)
+(x_2^*-x_3^*)\otimes(x_1^*-x_4^*)\wedge(x_4^*-x_7^*)
\\ & & {}
+(x_4^*-x_5^*)\otimes(x_3^*-x_6^*)\wedge(x_6^*-x_7^*)
+(x_1^*-x_6^*)\otimes(x_2^*-x_5^*)\wedge(x_5^*-x_7^*)
\\ & & {}
-x_7^*\otimes(x_2^*-x_4^*)\wedge(x_4^*-x_6^*)
+(x_4^*-x_5^*)\otimes x_1^*\wedge x_2^*
\\ & & {}
-(x_1^*-x_6^*)\otimes x_3^*\wedge x_4^*
+(x_2^*-x_3^*)\otimes x_5^*\wedge x_6^*
,\\ T_1 & = &
x_7^*\otimes(x_1^*-x_3^*)\wedge(x_3^*-x_5^*)
-x_3^*\otimes(x_1^*-x_4^*)\wedge(x_4^*-x_7^*)
\\ & & {}
-x_5^*\otimes(x_3^*-x_6^*)\wedge(x_6^*-x_7^*)
+x_1^*\otimes(x_2^*-x_5^*)\wedge(x_5^*-x_7^*)
\\ & & {}
-(x_5^*-x_7^*)\otimes x_1^*\wedge x_2^*
-(x_1^*-x_7^*)\otimes x_3^*\wedge x_4^*
-(x_3^*-x_7^*)\otimes x_5^*\wedge x_6^*
,\\ T_2 & = &
x_2^*\otimes(x_1^*-x_3^*)\wedge(x_3^*-x_5^*)
+(x_2^*-x_5^*)\otimes(x_3^*-x_6^*)\wedge(x_6^*-x_7^*)
\\ & & {}
+(x_3^*-x_6^*)\otimes(x_2^*-x_5^*)\wedge(x_5^*-x_7^*)
-x_3^*\otimes(x_2^*-x_4^*)\wedge(x_4^*-x_6^*)
\\ & & {}
+(x_3^*-x_5^*)\otimes x_1^*\wedge x_2^*
-(x_2^*-x_6^*)\otimes x_3^*\wedge x_4^*
+(x_2^*-x_3^*)\otimes x_5^*\wedge x_6^*
,\\ T_3 & = &
x_6^*\otimes(x_1^*-x_3^*)\wedge(x_3^*-x_5^*)
-(x_3^*-x_6^*)\otimes(x_1^*-x_4^*)\wedge(x_4^*-x_7^*)
\\ & & {}
-(x_1^*-x_4^*)\otimes(x_3^*-x_6^*)\wedge(x_6^*-x_7^*)
-x_1^*\otimes(x_2^*-x_4^*)\wedge(x_4^*-x_6^*)
\\ & & {}
+(x_4^*-x_6^*)\otimes x_1^*\wedge x_2^*
-(x_1^*-x_6^*)\otimes x_3^*\wedge x_4^*
+(x_1^*-x_3^*)\otimes x_5^*\wedge x_6^*
,\\ T_4 & = &
x_4^*\otimes(x_1^*-x_3^*)\wedge(x_3^*-x_5^*)
+(x_2^*-x_5^*)\otimes(x_1^*-x_4^*)\wedge(x_4^*-x_7^*)
\\ & & {}
+(x_1^*-x_4^*)\otimes(x_2^*-x_5^*)\wedge(x_5^*-x_7^*)
-x_5^*\otimes(x_2^*-x_4^*)\wedge(x_4^*-x_6^*)
\\ & & {}
+(x_4^*-x_5^*)\otimes x_1^*\wedge x_2^*
-(x_1^*-x_5^*)\otimes x_3^*\wedge x_4^*
+(x_2^*-x_4^*)\otimes x_5^*\wedge x_6^* .
\end{eqnarray*}
Moreover, we have $(R_3^\perp)^\perp=R_3$; hence $P_3$ is a free
Abelian group and $P_3^*\cong R_3^\perp$.
Therefore, it is enough to show that $U^\perp=\im\tilde{\tau}^*$.

We identify $(\Hom(R_2,P_3))^*$ with $R_2\otimes P_3^*$.
It is easy to check that
\begin{gather*}
I(p_{135})  =  \tilde{\tau}^*(r(1,p_{135})\otimes S_{135}),\quad
I(p_{147})  =  \tilde{\tau}^*(r(1,p_{147})\otimes S_{147}),\\
J(p_{16})  =  \tilde{\tau}^*(r(1,p_{16})\otimes T_3),\quad
K_1(p_{135})  =  \tilde{\tau}^*(r(1,p_{135})\otimes(T_0-T_3)),\\
K_2(p_{135})  =  \tilde{\tau}^*(r(1,p_{135})\otimes T_3),\quad
K_1(p_{147})  =  -\tilde{\tau}^*(r(1,p_{147})\otimes T_0),\\
K_2(p_{147})  =  -\tilde{\tau}^*(r(1,p_{147})\otimes T_3).
\end{gather*}
Hence $(U^\perp)_p\subset\im\tilde{\tau}^*$ for $p=p_{135},p_{147},p_{16}$.
Now we can prove that it is true for any $p\in\ppp_0$ by using the fact that the group generated by the permutations $(16)(25)(34)$ and $(135)(246)$ naturally acts on $C_8$ by automorphisms. Hence
$U^\perp\subset\im\tilde{\tau}^*$. Combining this fact with Proposition~\ref{inv_ker_tau}, we complete the proof of the lemma. }

\begin{lemma} The equality $B=\tilde{\tau}^{-1}(\bar{\delta}\Hom(H,P_2))$ holds true for the incidence structure $C_8$.
\end{lemma}

\prf{
We first show that if $\bar{\delta}f\in\tilde{\tau}A$ for some
$f\in\Hom(H,P_2)$, then

(1)~$(\omega_{ij},f(x_k))=0$ for all
$i,j,k$ such that $l_i\cap l_j\prec l_0$, $k\neq i,j$, and

(2)~$(\omega_{ijk},f(x_m))=0$ for all $i,j,k,m$ such that
$l_i\succ l_j\cap l_k\in\ppp_0$, $m\neq i,j,k$.

Assume that (1) doesn't hold; then, since no four pairwise distinct lines of
$C_8$ meet at one point, we have $p=l_i\cap l_k\in\ppp_0$ and
$(S_{ij},\bar{\delta}f(r(k,p)))=(\omega_{ij},f(x_k))\neq0$.
But $(S_{ij},\tilde{\tau}a(r(k,p)))=0$ for any $a\in A$, and
hence (1) must hold.

Assume that (2) doesn't hold; then we have either $p=l_i\cap l_m\in\ppp_0$
or $p'=l_j\cap l_m\in\ppp_0$. Due to the order symmetry of $(ijk)$, we
may suppose that $p=l_i\cap l_k\in\ppp_0$. Then
$(S_{ijk},\bar{\delta}f(r(m,p)))=(\omega_{ijk},f(x_m))\neq0$.
But $(S_{ijk},\tilde{\tau}a(r(m,p)))=0$ for any $a\in A$, and
hence (2) must hold as well.

It follows that we can choose $\tilde{f}\in\Hom(H,\Lambda^2H)$ with
$f(x)=\tilde{f}(x)+R_2$ in such a way that $\tilde{x_k}=x_k\wedge
b(k)$ for some $b:\{1,\dots,7\}\to H$. Hence
$\bar{\delta}f=\tilde{\tau}a$, where $a\in B$, $a(k,p)=b(k)$
for any $k$ and $p$. }

It follows from the last two lemmas that the homomorphism $\bar{\tau}:W\to T$
is injective; hence $\tau a_0\neq0$. This concludes the proof of Theorem \ref{inv_mclane}. }

In the course of computations, we produced a nonzero homomorphic image of the distinguishing invariant $\tau a_0$ in $\zz/3\zz$. Some additional computations show that the element $\tau a_0$ itself is of order $3$ in the group $T$ and thus belongs to its torsion subgroup. This is no accident: the invariant must vanish after tensoring with~$\qq$, since the Malcev completions of the fundamental groups are canonically isomorphic (see~\cite{Ko}).

Note that the groups $\po(X^+)$ and $\po(X^-)$ are, in fact, isomorphic, because complex conjugation induces a homeomorphism $X^+\to X^-$. But since this homeomorphism reverses the coorientation of lines,
the corresponding isomorphism of the first homology groups maps the canonical
basis of one of them to the negative canonical basis of the other.

Let us try to describe the group automorphisms of $G=\po(X^+)\cong\po(X^-)$.

As mentioned above, $\Aut(C_8)\cong GL_2(\ff_3)$. It is readily seen that the automorphisms corresponding to the elements of $SL_2(\ff_3)$ transform the realization $\phi^{\pm}$ into projectively equivalent ones, while the remaining automorphisms transform $\phi^{+}$ into realizations projectively equivalent to $\phi^{-}$, and vice versa.

Let $g\in GL_2(\ff_3)$. By $\sigma_g$ we denote a permutation of $\{0,1,\dots,7\}$ such that each line $l_i$ is mapped to $l_{\sigma_g(i)}$ under the automorphism of $C_8$ corresponding to $g$. We also choose a loop $w_0\in G$ passing around the line at infinity $l_0$ in the positive direction. Since any projective equivalence determines an isomorphism of fundamental groups, for each $g\in SL_2(\ff_3)$, we obtain an automorphism
$\eta_g$ of~$G$ which maps each element $w_i\in G$ ($i=0,1\dots,7$) to a loop passing around $l_{\sigma_g(i)}$ in the positive direction, that is, to an element conjugate to $w_{\sigma_g(i)}$ in
$G$. Similarly, by using the isomorphism of $\po(X^+)$ and $\po(X^-)$ induced by complex conjugation, for each $g\in GL_2(\ff_3)\setminus SL_2(\ff_3)$, we obtain an automorphism $\eta_g$ of~$G$ which maps each element $w_i\in G$ ($i=0,1\dots,7$) to a loop passing around $l_{\sigma_g(i)}$ in the negative direction, that is, to an element conjugate to $w_{\sigma_g(i)}^{-1}$ in $G$.

Let $H=G/\gamma_2G$. Since automorphisms of $G$ map $\gamma_2G$ onto itself, we can reduce any automorphism $\eta\in\Aut G$ modulo $\gamma_2G$ so as to obtain $\bar\eta\in\Aut H$. Let $W=\{\bar\eta\mid\eta\in\Aut G\}$.

\begin{theorem}\label{aut_mclane}
The correspondence $g\mapsto\bar\eta_g$ is a group isomorphism $GL_2(\ff_3)\to W$.
\end{theorem}

\begin{proof}
Consider the element $x_0=-\sum_{i=1}^7x_i$ in $H$. It is readily seen that
$$
\bar\eta_g(x_i)=
  \begin{cases}
    x_{\sigma_g(i)} & \text{if $g\in
    SL_2(\ff_3)$}, \\
    -x_{\sigma_g(i)} & \text{if $g\notin
    SL_2(\ff_3)$}
  \end{cases}
$$
for $i=0,1,\dots,7$. Hence we see that the correspondence $g\mapsto\bar\eta_g$
is an injective homomorphism. Let us prove that it is surjective.

Take any automorphism $\eta\in\Aut G$. Let us extend $\bar\eta\in\Aut H$ to an automorphism of
$\Lambda^2H$. Clearly, this automorphism preserves the subgroup $R_2=\Ker(\Lambda^2H\to P_2)$, where $P_2=\gamma_2G/\gamma_3G$ (see Section~2). Consider $R_2^\perp\subseteq\Lambda^2H^*$. As we know, the Abelian group $R_2^\perp$ is generated by the elements $\omega_{ijk}=(x_i^*-x_j^*)\wedge(x_j^*-x_k^*)
=x_i^*\wedge x_j^*+x_j^*\wedge x_k^*+ x_k^*\wedge x_i^*$, where $l_i\succ l_j\cap l_k\in\ppp_0$, and $\omega_{0ij}=\omega_{ij}=x_i^*\wedge x_j^*$, where $l_0\succ l_i\cap l_j$.

\begin{lemma}
In the case of the incidence structure $C_8$, each decomposable element of $R_2^\perp\subseteq\Lambda^2H^*$
is proportional to one of the $\omega_{ijk}$, where $p_{ijk}\in\ppp$.
\end{lemma}

\begin{proof}
Since $H^*$ is a free Abelian group, we can perform computations over $\qq$ in $H^*\otimes\qq$. Note that the elements $\omega_{ijk}$ are linearly independent, and the group $\Aut(C_8)$ maps the one-dimensional subspaces of $H^*\otimes\qq$ generated by these elements to each other. Moreover, the action of this group on the set $\{\left<\omega_{ijk}\right>\mid p_{ijk}\in\ppp\}$ of these subspaces is transitive. Thus, it suffices to show that if a decomposable element $\omega$ is expressed as a linear combination of $\omega_{ijk}$ in which the coefficient of $\omega_{012}$ is $1$, then $\omega=\omega_{012}$.

Indeed, let us represent $\omega$ in the form $\omega=\sum_{0<i<j}a_{ij}x_i^*\wedge x_j^*$. Then $a_{12}=1$, and $\omega$ can be decomposed as $\omega=(x_1^*+\sum_{i=3}^7c_ix_i^*)\wedge(x_2^*+\sum_{i=3}^7d_ix_i^*)$. On the other hand, $\omega$ is a linear combination of the forms $\omega_{ijk}$. Comparing these two decompositions, we readily see that $\omega=\omega_{012}$.
\end{proof}

Let $V=H\otimes\qq$. Since $H$ is a free $\zz$-module, we can regard it as a submodule of $V$. The vector space $V$ can be treated as $\qq^8/\left<e_0+\dots+e_7\right>$, where $(e_0,\dots,e_7)$ is the standard basis of $\qq^8$. Here $x_i=e_i+\left<e_0+\dots+e_7\right>$ for all
$i=0,\dots,7$. Clearly, the action of the automorphism group of the incidence structure
$C_8$ on $V$ comes from its action by the permutations $\sigma_g$ of the basis vectors in~$\qq^8$.

The dual space $V^*$ can be represented as the subspace in $(\qq^8)^*$ generated by $e_i^*-e_j^*$. Its basis dual to the basis
$(x_1,\dots,x_7)$ of $V$ consists of the elements $x_i^*=e_i^*-e_0^*$ ($i=1,\dots,7$). As we have shown, the automorphism $\bar\eta$ transforms each form $\omega_{ijk}=(e_i^*-e_j^*)\wedge(e_j^*-e_k^*)$, where $p_{ijk}\in\ppp$, into a form of the same type (up to sign). Consider the kernel of this form in $V$. It consists of all vectors in which the coefficients of $e_i$, $e_j$, and $e_k$ are equal. It is easy to check that minimal nonzero subspaces in $V$ which can be expressed as the intersections of the kernels of some of the $\omega_{ijk}$ are precisely the one-dimensional subspaces generated by the $x_s$ ($s=0,1,\dots,7$). Therefore, the automorphism $\bar\eta$ maps each $x_s$ into $\pm x_{\sigma(s)}$ for some permutation $\sigma\in S_8$. Since the vectors $x_s$ satisfy the relation $\sum x_s=0$, all of the signs ``$\pm$'' in front of $x_{\sigma(s)}$ are the same.
Further, representing each $\left<x_s\right>$ as the intersection of the kernels of some forms $\omega_{ijk}$,
it is easy to show that $\sigma=\sigma_g$ for some $g\in GL_2(\ff_3)$, and the sign is ``$+$'' for $g\in SL_2(\ff_3)$ and ``$-$'' otherwise. The latter assertion follows from the nonexistence of an automorphism 
$\eta\colon G\to G$ such that $\bar\eta(x_i)=-x_i$ for all $i=0,1,\dots,7$
(otherwise, applying complex conjugation, we would obtain a counterexample to Theorem~\ref{inv_mclane}).
\end{proof}

\begin{remark}
Theorem~\ref{aut_mclane} essentially means that in the case of the MacLane arrangement the standard basis in the first homology group, which is the Abelianization of the fundamental group, is ``almost'' an invariant of the fundamental group itself. To be more precise, the set of cyclic subgroups generated by the elements of the standard basis and by $x_0$ (the cycle around the line at infinity) is invariant under automorphisms. Moreover, each automorphism maps the set $\{x_0,\dots,x_7\}$ either onto itself or onto set $\{-x_0,\dots,-x_7\}$. This ``structural rigidity'' of the standard basis for the MacLane arrangement will be used in the next section.

Not all line arrangements possess this property. For example, in the case of lines in general position,
the fundamental group is Abelian, and the standard basis of the first homology group cannot be reconstructed in any sense from the fundamental group.\footnote{I am grateful to the referee, who advised me to add this remark.}
\end{remark}

\section{Combinatorially Equivalent Arrangements with Different Fundamental Groups}

We define an incidence structure $C_{13}$ by gluing together two MacLane
structures in the following way. Let $C'_8=(\clll',\ppp',\succ)$
and $C''_8=(\clll'',\ppp'',\succ)$ be two copies of the MacLane structure
($\clll'=\{l_0',\dots,l_7'\}$, $\ppp'=\{p_{012}',\dots,p_{45}'\}$,
$\clll''=\{l''_0,\dots,l''_7\}$, $\ppp''=\{p''_{012},\dots,p''_{45}\}$).
Let us introduce an equivalence relation on $\clll'\cup\clll''$
by setting $l_0'\sim l''_0$, $l_1'\sim l''_1$, and $l_2'\sim l''_2$,
and an equivalence relation on $\ppp\cup\ppp'$
by setting $p_{012}'\sim p''_{012}$. Let $\ppp'''=\{p'''_{ij}\mid
i,j\in\{3,\dots,7\}\}$. We set $\hat{\clll}=(\clll'\cup\clll'')/{\sim}$,
$\hat{\ppp}=((\ppp'\cup\ppp'')/{\sim})\cup\ppp'''$. The incidence relation
between $\hat{\clll}$ and $(\ppp\cup\ppp')/{\sim}$ is induced by those of
$C'_8$, $C''_8$, and the incidence relation between $\hat{\clll}$ and
$\ppp'''$ is $p'''_{ij}\prec l_i',l''_j$. Let $C_{13}=(\hat{\clll},\hat{\ppp},\succ)$.
We denote the equivalence class of $l'_i$ ($i=0,1,2$) by $\hat{l_i}$ and
the equivalence class of $p'_{012}$ by $\hat{p}_{012}$.

Let $\psi$ be a projective transformation of $\cc P^2$
such that $\psi\phi^+l_i=\phi^+l_i(=\phi^-l_i)$
for $i=0,1,2$. This condition means that $\psi$ acts as $(z_0:z_1:z_2)\mapsto(z_0:z_1:az_0+bz_1+cz_2)$, where
$a,b,c\in\cc, c\ne0$. Assume that none of the lines $\psi\phi^\pm l_i$,
$i=3,\dots,7$, passes through any of the intersection points of the lines
$\phi^+l_j$, $j=0,1,\dots,7$. In this case we say that $\psi$ is generic.
Let us construct realizations $\phi^{++}$ and $\phi^{+-}$ of $C_{13}$
by gluing together the realizations $\phi^+$ and $\phi^-$ of $C_8$ as follows:
we set $\phi^{+\pm}\hat{l_i}=\phi^+l_i$ for $i=0,1,2$,
$\phi^{+\pm}l'_i=\phi^+l_i$ and $\phi^{+\pm}l''_i=\psi\phi^\pm l_i$
for $i=3,\dots,7$. For generic $\psi$, the intersection points of
complex lines $\phi^{+\pm}\hat l$ (where $\hat l\in\hat{\clll}$)
are in one-to-one correspondence with the elements of $\hat{\ppp}$,
hence we obtain realizations of $C_{13}$. Let $X^{+\pm}_\psi=\cc
P^2\setminus\bigcup\phi^{+\pm}_\psi l$. Note that the space of generic transformations $\psi$
is connected, and hence the fundamental groups $G^{++}=\po(X^{++}_\psi)$ and
$G^{+-}=\po(X^{+-}_\psi)$ are independent of $\psi$ up to an isomorphism preserving the canonical basis in $G^{+\pm}/\g2G^{+\pm}$. Therefore, we can omit the index~$\psi$.

\begin{theorem}
The groups $G^{++}$ and $G^{+-}$ are not isomorphic.
\end{theorem}

\prf{
Let $\hat G$ be one of the groups $G^{++}$ and $G^{+-}$, and let $\hat
x_1,\hat x_2$, $x'_3,\dots,x'_7$, $x_3'',\dots,x_7''$ denote the images of the standard generators of $\hat G$ in $\hat H = \hat G/\gamma_2\hat G$. We set
$\hat x_0=-\hat x_1-\hat x_2-x'_3-\dots-x'_7-x_3''-\dots-x_7''$.

\begin{lemma}\label{mclane_pr_inc}
As in the previous section, let $G=\po(X^+)$. There exist concordant projections
$\pi'\colon\hat G\to G$, $\pi''\colon\hat G\to G$ and injections $\nu'\colon
G\to\hat G$, $\nu''\colon G\to\hat G$ such that the following relations hold modulo the second term of the lower central series:
\begin{eqnarray*}
   &\bar\pi'(\hat x_i)= x_i\quad (0\les i\les 2), \quad \bar\pi'(x'_i) = x_i\text{ \ and \ }
   \bar\pi'(x''_i) = 0 \quad (3\les i\les 7), \\
   &\bar\nu'(x_i) = \hat x_i\quad (0\les i\les 2), \quad \bar\nu'(x_i) = x'_i \quad (3\les i\les 7),
\end{eqnarray*}
as well as
\begin{eqnarray*}
   &\bar\pi''(\hat x_i)= x_i\quad (0\les i\les 2), \quad \bar\pi''(x'_i) = 0\text{ \ and \ }
   \bar\pi''(x''_i) = x_i \quad (3\les i\les 7), \\
   &\bar\nu''(x_i) = \hat x_i\quad (0\les i\les 2), \quad \bar\nu''(x_i) = x''_i \quad (3\les i\les 7),
\end{eqnarray*}
for $\hat G=G^{++}$, and
\begin{eqnarray*}
   &\bar\pi''(\hat x_i)= -x_i\quad (0\les i\les 2), \quad \bar\pi''(x'_i) = 0\text{ \ and \ }
   \bar\pi''(x''_i) = -x_i \quad (3\les i\les 7), \\
   &\bar\nu''(x_i) = -\hat x_i\quad (0\les i\les 2), \quad \bar\nu''(x_i) = -x''_i \quad (3\les i\les 7),
\end{eqnarray*}
for $\hat G=G^{+-}$.
\end{lemma}

\prf{
The required projections are defined in an obvious way, by removing the lines $\psi\phi^\pm l_i$ for
$i=3,\dots,7$ (or the lines $\phi^+ l_i$ for $i=3,\dots,7$). To construct the injections, we choose a generic projective transformation $\psi$ so that the lines $\psi\phi^\pm l_i$ for $i=3,\dots,7$ do not intersect some open ball containing all intersection points of the affine lines $\phi^+ l_i$ for $i=1,\dots,7$, and the lines $\phi^+ l_i$ for $i=3,\dots,7$ do not intersect some open ball containing all intersection points of affine lines $\psi\phi^\pm l_i$ for $i=1,\dots,7$. It is readily seen that the complements of lines $\phi^+ l_i$ and $\psi\phi^\pm l_i$ in these open balls are homeomorphic to $X^+$ and $X^\pm$, respectively.
}

\begin{lemma}\label{aut_mclane_pp_pm}
Let $\eta\colon\hat G\to\hat G$ be an automorphism, and let $\bar\eta\colon\hat
H\to\hat H$ be its reduction modulo $\gamma_2\hat G$. Then there exists an element $g\in GL_2(\ff_3)$, $\sigma_g(\{0,1,2\})=\{0,1,2\}$, and a number $\ep=\pm1$ such that $\bar\eta(\hat x_i)=\ep\hat x_{\sigma_g(i)}$ for $0\les i\les2$ and either $\bar\eta(x'_i)=\ep x'_{\sigma_g(i)}$ and $\bar\eta(x''_i)=\ep x''_{\sigma_g(i)}$ or $\bar\eta(x'_i)=\ep x''_{\sigma_g(i)}$ and
$\bar\eta(x'_i)=\ep x''_{\sigma_g(i)}$ for $3\les i\les7$. The same is true for any isomorphism $G^{++}\to G^{+-}$.
\end{lemma}

\prf{
The required assertion follows from Theorem~\ref{aut_mclane} and Lemma~\ref{mclane_pr_inc}.
}

\begin{lemma}\label{aut_mclane_pp}
There exists an automorphism $\eta\colon G^{++}\to G^{++}$ such that
$\bar\eta(\hat x_i)=\hat x_i$ for $0\les i\les2$ and $\bar\eta(x'_i)=x''_i$,
$\bar\eta(x''_i)=x'_i$ for $3\les i\les7$.
\end{lemma}

\prf{
For any generic projective transformation $\psi$ of the required form, the inverse transformation
$\psi^{-1}$ is also generic; therefore, $\psi$ induces an isomorphism $\nu\colon\po(X^{++}_{\psi^{-1}})\to\po(X^{++}_\psi)$. Consider  an isomorphism $\mu\colon\po(X^{++}_\psi)\to\po(X^{++}_{\psi^{-1}})$ defined by homotopy along some path connecting $\psi$ with $\psi^{-1}$ in the space of generic transformations. Since $\mu$ is concordant with the canonical bases, while $\nu$ interchanges the images of $l'_i$ and
$l''_i$, their composition $\eta=\nu\mu\colon\po(X^{++}_\psi)\to\po(X^{++}_\psi)$
is an automorphism with the required properties.
}

\begin{lemma}\label{aut_mclane_pm}
There exists no automorphism $\eta\colon G^{+-}\to G^{+-}$ such that
$\bar\eta(\hat x_i)=\hat x_i$ for $0\les i\les2$ and $\bar\eta(x'_i)=x''_i$,
$\bar\eta(x''_i)=x'_i$ for $3\les i\les7$.
\end{lemma}

\prf{
The required assertion follows from Theorem~\ref{aut_mclane} and Lemma~\ref{mclane_pr_inc}.
}

The assertion of the theorem readily follows from Lemmas~\ref{aut_mclane_pp_pm},~\ref{aut_mclane_pp}, and~\ref{aut_mclane_pm}.
}

\end{document}